\documentclass[a4paper]{amsart}
\usepackage[T1]{fontenc}
\usepackage[latin1]{inputenc}

\makeatletter

\usepackage[T1]{fontenc}
\usepackage[latin1]{inputenc}
\usepackage{amssymb}
\usepackage[T1]{fontenc}
\usepackage{graphics}
\usepackage{verbatim}
\makeatletter
\usepackage{graphicx}
\usepackage{amsthm}
\usepackage{amssymb}
\usepackage{amsfonts}
\usepackage{amsmath}
\usepackage{amsfonts}

\newcommand{\R}{\mathbb{R}}

\newcommand{\K}{\mathcal{K}}
\newtheorem{axiom}{Axiom}

\newtheorem{corollary}[axiom]{Corollary}
\newtheorem{definition}[axiom]{Definition}

\newtheorem{lemma}[axiom]{Lemma}

\newtheorem{proposition}[axiom]{Proposition}
\newtheorem{theorem}[axiom]{Theorem}
\theoremstyle{definition}
\newtheorem{remark}[axiom]{Remark}
\makeatother
\begin{document}
\title{On the shrinkage behavior of Partial Least Squares Regression}
\author{Nicole Kr\"amer}

\maketitle
\begin{abstract}
We present a formula for the shrinkage factors of the Partial Least Squares
regression  estimator and
deduce some of their properties, in particular the known fact that some of
the factors are $>1$. We investigate the effect of shrinkage factors for the
Mean Squared error of linear estimators and illustrate that we cannot extend
the results to nonlinear estimators. In particular, shrinkage factors $>1$ do not
automatically lead to a poorer Mean Squared Error. We investigate empirically the
effect of bounding the the absolute value of the Partial Least Squares shrinkage factors by $1$. 
\end{abstract}

\keywords{\textbf{Keywords}  Partial Least Squares, shrinkage estimators}\\

\subjclass{\textbf{AMS classification} 62J07, G2H99}
\section{Introduction}
We investigate the shrinkage properties of the Partial Least
Squares (PLS) regression estimator. It is known (e.g.
\cite{Frank9301}) that we can express the PLS estimator obtained
after $m$ steps in the following way:
\begin{align*}
\hat \beta_{PLS} ^{(m)}=\sum_{i=1} ^p f^{(m)}(\lambda_i) \cdot  z_i \,,
\end{align*}
where $z_i$ is the component of the Ordinary Least Squares (OLS) estimator along the $i$th
principal component of the covariance  matrix $X^t X$ and $\lambda_i$ is the
corresponding
eigenvalue. The quantities $f^{(m)}(\lambda_i)$ are called shrinkage factors. We
show that these factors are determined by a tridiagonal matrix (which
depends on the input--output matrix $(X,y)$) and can be
calculated in a recursive way. Combining the results of \cite{Butler0001}
and \cite{Phatak0301},  we give a simpler and clearer proof of the shape of the
shrinkage factors of PLS and derive some of their properties. In particular, we show that
some of the values $f^{(m)}(\lambda_i)$ are greater than $1$ (this was first
proved in \cite{Butler0001}). \\

We argue that these "peculiar
shrinkage properties" \cite{Butler0001} do not necessarily imply that the Mean Squared Error
(MSE) of the PLS
estimator is worse compared to the MSE of the OLS estimator:
In the case of deterministic shrinkage factors, i.e. factors that do not
depend on the output $y$, any value
$\left |f^{(m)}\left(\lambda_i\right)\right|>1$ is of course undesirable. But in
the case of PLS, the shrinkage factors are stochastic -- they also depend on $y\,$. Even if $P\left( \left
|f^{(m)}\left(\lambda_i\right)\right|>1\right)=1$ we cannot conclude that the MSE
is worse than the MSE of the OLS estimator. In particular, bounding the
absolute value of the shrinkage factor by $1$ does not
automatically yield a lower MSE, in disagreement to what was
conjectured in e.g. \cite{Frank9301}.\\

Having issued this warning, we explore whether bounding the shrinkage
factors leads to a lower MSE or not. It is very difficult to
derive theoretical results, as the quantities of interest -
$\hat\beta_{PLS} ^{(m)}$ and $f^{(m)}(\lambda_i)$ respectively - depend on $y$ in a
complicated, 
nonlinear way. As a substitute, we study the problem on several artificial data sets
and one real
world example. It turns out that in most cases the MSE of
the bounded version of PLS is indeed smaller than the one of PLS, although
the improvement is tiny. \\

The paper is organized as follows: In section \ref{pre} we
introduce the notation and in section \ref{kspaces} we recall some propertie
of Krylov spaces. In section \ref{pls} we define the PLS
estimator and in section \ref{sectiontri} we provide mathematical results that are needed in the
rest of the paper. After explaining the notion of shrinkage in
section \ref{shrinkage} we derive the formulas for the PLS
shrinkage factors in section \ref{shrinkpls} and derive some of their
properties. In section \ref{Simulation}, we report on the results of the
experiments. The paper ends with a conclusion.
\section{Preliminaries}
\label{pre}
We consider the multivariate linear regression model
\begin{eqnarray}
\label{linreg}
y&=&X\beta +\varepsilon\,
\end{eqnarray}
with
\begin{eqnarray*}
\text{Cov}\left(y\right)&=& \sigma^2 \cdot \text{Id}\,.
\end{eqnarray*}
The numbers of variables is $p$, the number of examples is $n\,$. For
simplicity, we assume that $X$ and $y$ are scaled to have zero mean, so we
do not have to worry about intercepts. We have
\begin{eqnarray*}
X &\in& \R^{n\times p}\,,\\
A:=X^tX &\in& \R^{p\times p} \,,\\
y &\in& \R^n\,,\\
b:=X^ty &\in& \R^p\,.\\
\end{eqnarray*}
We set $p^*=\text{rk}\left(A\right)=\text{rk}\left(X\right)$. The singular value decomposition of $X$ is of the form
\begin{eqnarray*}
X&=& V\Sigma U^t\,
\end{eqnarray*}
with
\begin{eqnarray*}
V&\in &\mathbb{R}^{n\times p}\,\\
\Sigma=\text{diag}\left(\sqrt{\lambda_1},\ldots,\sqrt{\lambda_p}\right)&\in& \mathbb{R}^{p\times p}
\,,\\
U&\in& \R^{p\times p} \,.\\
\end{eqnarray*}
We have $U^t  U=\text{Id}_p$ and $V^t  V=\text{Id}_p$. \\

Set $\Lambda=\Sigma^2\,$. The eigendecomposition of $A$ is 
\begin{eqnarray*}
A&=&U\Lambda U^t=\sum_{i=1} ^p \lambda_i u_i u_i ^t\,.
\end{eqnarray*}

The eigenvalues $\lambda_i$ of $A$ (and any other matrix) are  ordered in the following way:
\begin{eqnarray*}
\lambda_1 \geq \lambda_2\geq \ldots \geq \lambda_p \geq 0\,.
\end{eqnarray*}
The Moore-Penrose inverse of a matrix $M$ is denoted by $M^-$. \\

The Ordinary Least Squares (OLS) estimator $\hat \beta_{OLS}$ is the
solution of the optimization problem 
\begin{eqnarray*}
\text{arg} \min_{\beta} &\,& \|y-X\beta\|\,.
\end{eqnarray*} 
Set 
\begin{eqnarray}
\label{te}
t&=& \Sigma V^t y\,.
\end{eqnarray}
The OLS estimator is given by the formula
\begin{eqnarray*}
\hat \beta_{OLS}&=& \left(X^tX\right)^- X^t y\\
&=& U \Lambda^- U^t U\Sigma V^t y\\
&=& U\Lambda^- t\\
&=& \sum_{i=1} ^{p^*} \frac{v_i ^t y}{\sqrt{\lambda_i}} u_i\,.
\end{eqnarray*}
Set 
\begin{eqnarray*}
z_i&=& \frac{v_i ^t y}{\sqrt{\lambda_i}} u_i\,.
\end{eqnarray*}

Finally, we need a result on the shape of the Moore-Penrose inverse of a symmetric
matrix.
\begin{proposition}\label{penrose}
Let $B \in \mathbb{R}^{m\times m}$ be a symmetric matrix with
eigendecomposition 
\begin{eqnarray*}
B=S \Lambda S^t\,,
\end{eqnarray*}
with eigenvalues $\lambda_i$. Set
\begin{eqnarray*}
f_B(\lambda)&=&1-\prod_{\lambda_i\not=0} \left(1-\frac{\lambda}{\lambda_i}  \right)\,.
\end{eqnarray*}
As $f_B(0)=0$ we can write 
\begin{eqnarray*}
f_B(\lambda)&=&\lambda\cdot \pi_B(\lambda)\,. 
\end{eqnarray*}
Then
\begin{eqnarray*}
B^-&=&\pi_B(B)\,.
\end{eqnarray*}
\end{proposition}
\begin{proof}
The four properties that we have to check are
\begin{enumerate}
\item $\left( B B^-\right)^t=BB^-\,$,
\item  $\left(  B^- B\right)^t=B^- B\,$,
\item $BB^-B=B\,$,
\item $B^-BB^-=B^-\,$.
\end{enumerate}
As $B$ is symmetric, the polynomial $\pi_B(B)$ is symmetric as well, which
proves the first two conditions.
Next note that it suffices to prove the and  properties $3$ and $4$ for the diagonal matrix
\begin{eqnarray*}
\Lambda&=&\text{diag}\left(\lambda_1,\ldots,\lambda_k,0\ldots,0\right)
\end{eqnarray*}
 with $k=\text{rk}(B)$. This is true as 
\begin{eqnarray*}
B^-&=& \left(S \Lambda S^t\right)^-\\ 
&=& S \Lambda^- S^t\,.
\end{eqnarray*}
We have
\begin{eqnarray*}
\Lambda^-&=&\text{diag}\left(\lambda_1
^{-1},\ldots,\lambda_k ^{-1},0\ldots,0\right)\,.
\end{eqnarray*}
 The third property of the
Moore-Penrose inverse is $\Lambda= \Lambda \Lambda^- \Lambda$ which is
equivalent to $\lambda_i=\lambda_i \pi_B(\lambda_i) \lambda_i$ which is
obviously true. The fourth property follows as easily.
\end{proof}
\begin{remark}
The degree of the polynomial $\pi_B$ is $\text{rk}(B)-1\,$. The proposition
is valid no matter if we count the non-zero eigenvalues with or without
multiplicities. We count the eigenvalues with multiplicities in order to
connect the polynomial to the characteristical polynomial in the regular case: If $B$ is a regular
matrix, $\pi_B$ is linked to the characterictical polynomial $\chi_B$ in the
following way:
\begin{eqnarray*}
\lambda \cdot \pi_B (\lambda)&=&\frac{1}{\chi_B (0)} \chi_B(\lambda) +1\,. 
\end{eqnarray*}
\end{remark}
\section{Krylov spaces}
\label{kspaces}
Set 
\begin{eqnarray*}
K^{(m)}:=\left (A^0b,Ab,\ldots,A^{m-1}b\right ) &\in&  \R^{p\times m}\,.
\end{eqnarray*}
The columns of $K^{(m)}$ are called the Krylov sequence of $A$ and $b$. \\

The space spanned by the columns of $K^{(m)}$ is called the Krylov space of
$A$ and $b$ and denoted
by $\K^{(m)}$. We recall some basic facts on the dimension of the Krylov space that are
needed in the rest of the paper. Set
\begin{eqnarray*}
\mathcal{M}&:=\left\{\lambda_i|t_i\not=0\right\}
\end{eqnarray*}
( the vector $t$ is defined in (\ref{te})) and
\begin{eqnarray*}
m^*&:=&\left|\mathcal{M}\right|\,.
\end{eqnarray*}
\begin{lemma}
We have 
\begin{eqnarray*}
\dim \K^{\left(m^*\right)}&=&m^*\,.
\end{eqnarray*}
\end{lemma}
\begin{proof}
Suppose that 
\begin{eqnarray*}
\sum_{j=0} ^{m^*-1} \gamma_j A^j b &=&0
\end{eqnarray*}
for some $\gamma_0,\ldots,\gamma_{m^*-1} \in \mathbb{R}$. Using the
eigendecompostion of $A$ this equation is equivalent to 
\begin{eqnarray*}
U \left( \sum_{j=0} ^{m^*-1} \gamma_j \Lambda^j t \right)&=0
\end{eqnarray*}
As $U$ is an invertible matrix, this is equivalent to 
\begin{eqnarray*}
\sum_{j=0} ^{m^*-1} \gamma_j \lambda_i ^j t_i&=&0
\end{eqnarray*} 
for $i=1,\ldots,p$. Hence, each element $\lambda_i \in \mathcal{M}$ is a
zero of the polynomial
\begin{eqnarray*}
\sum_{j=0} ^{m^*-1} \gamma_j \lambda^j\,.
\end{eqnarray*}
This is  a polynomial of degree $\leq
m^*-1\,$. as it has $m^*=|\mathcal{M}|$ different zeroes, it must be
trivial, i.e. $\gamma_j=0$.
\end{proof}
\begin{lemma}
If $m>m^*$ we have  $\dim \mathcal{K}^{(m)}=m^*\,$.
\end{lemma}
\begin{proof}
It is clear that $ \dim \mathcal{K}^{(m)}\geq m^*\,$ as
$\mathcal{K}^{(m^*)}\subset \mathcal{K}^{(m)}\,$. Assume that there is a set $S$
of $m^*+1$ linear independent vectors in the Krylov sequence $K^{(m)}$. Set 
\begin{eqnarray*}
I&=&\{i\in \{1,\ldots m\}|A^{i-1}b  \in S\}\,.
\end{eqnarray*}
 Hence $|I|=m^*+1 $. The condition
that $S$ is linear independent is  equivalent to the following: There is
no nontrivial polynomial
\begin{eqnarray*}
g(\lambda)&=&\sum_{i\in I}  \gamma_i \lambda^i
\end{eqnarray*}
such that  
\begin{eqnarray}\label{dep}
g\left(\lambda_i\right)&=&0
\end{eqnarray}
 for $\lambda_i \in \mathcal{M}$. As
the polynomial $g$ is of degree $|I|=m^* +1$ and $|\mathcal{M}|=m^*$, there
is always a nontrivial solution of equation (\ref{dep}).
\end{proof}
We sum up the two results: 
\begin{proposition}
\label{dimkm}
We have
\begin{eqnarray*}
\dim \mathcal{K}^{(m)}&=&\begin{cases}
m & m\leq m^* \\
m^* & m>m^*
\end{cases}\,.
\end{eqnarray*}
In particular
\begin{eqnarray}
\label{chainkm}
\dim \K^{(m^*)}= \dim \K^{(m*+1)}=\ldots=\dim \K^{(p)}=m^*\,.
\end{eqnarray}
\end{proposition}

\section{Partial Least Squares}\label{pls}
It is not our aim to give an introduction to the Partial Least
Squares (PLS) method and refer to \cite{Hoeskuldsson8801}. We take a
purely algebraic  point of view as in \cite{Helland8801}. The
$PLS$ estimator $\hat \beta^{(m)} _{PLS}$ is the solution of the constrained minimization
problem
\begin{eqnarray*}
\text{arg}\min_{\beta} & \|y-X\beta\| \\
\text{s.t.}& \beta \in \K^{(m)} \,.
\end{eqnarray*}
We call $m$ the number of steps of PLS. It follows that any solution of this problem is of the form $\hat{\beta}=K^{(m)}
\hat z$ where $\hat z$ is the solution of the unconstrained problem
\begin{eqnarray*}
\text{arg}\min_{z}& \|y-XK^{(m)} z\|\,.
\end{eqnarray*}

Plugging this into the formula for the OLS estimator  (cf.  
section \ref{pre}) we get
\begin{proposition}[\cite{Helland8801}]
The PLS estimator obtained after $m$ steps can be expressed in the following way:
\begin{eqnarray}
\label{bhat}
\hat \beta_{PLS} ^{(m)} &= &K^{(m)} \left[\left(K^{(m)}\right) ^t A K^{(m)} \right]^{-}  \left(K^{(m)}\right)^t b\,.
\end{eqnarray}
\end{proposition}
It should be clear that we can replace the matrix $K^{(m)}$ in equation
(\ref{bhat}) by any matrix $W^{(m)}$, as long as its columns span the space
$\K^{(m)}$. In fact, in the NIPALS algorithm (see \cite{Helland8801}), an orthogonal basis of $\K^{(m)}$ is calculated with the help
of the Gram-Schmidt procedure. Denote by
\begin{eqnarray}
\label{gram}
W^{(m)}&=&\left(w_1,\ldots,w_m\right)
\end{eqnarray}
this orthogonal basis of $\K^{(m)}$. Of course, this basis only exists if $dim(\K^{(m)})=m$, which might
not be true for all $m\leq p$. The maximal number  for which this holds
is $m^*\,$ (see proposition \ref{dimkm}). Note however that 
\begin{eqnarray*}
\K^{(m^*-1)}\subset\K^{(m^*}=\K^{(m*+1)}=\ldots=\K^{(p)}
\end{eqnarray*}
(see  (\ref{chainkm})) and the solution of the optimization problem
does not change anymore. Hence for the rest of the paper, we make the assumption that
\begin{eqnarray}
\label{ass}
\dim \K^{(m)}=m\,.
\end{eqnarray}
\begin{remark}
We have
\begin{eqnarray*}
\hat \beta^{(m^*)} _{PLS}&=& \hat \beta_{OLS}\,. 
\end{eqnarray*}
\end{remark}
\begin{proof}
We show that $\hat \beta_{OLS} \in \mathcal{K}^{(m^*)}$. By definition
\begin{eqnarray*}
\hat \beta_{OLS}&=& U \Lambda^- t\\
&\stackrel{\ref{penrose}}{=}& U \pi_{\Lambda} (\Lambda) t\,.
\end{eqnarray*}
with $\text{deg} \pi_{\Lambda}= p^*-1\,$ (recall that $p^*$ is the rank of $A$). On the other hand, any vector
$v\in \mathbb{R}^{p^*}$  lies in  $\K^{(p^*)}$ if and only if there is a polynomial $g$ of degree $\leq p^*-1$ such that
\begin{eqnarray*}
v&=& g(A)b\\
&=&g\left(U\Lambda U^t \right) Ut\\
&=& U g(\Lambda) t\,.
\end{eqnarray*} 
It follows that  $\hat \beta_{OLS} \in \K^{(p^*)}$. As $p^* \geq m^*$ we
have $\K^{(p^*)}= \K^{(m^*)}\,$.  
\end{proof}
Set
\begin{eqnarray*}
T^{(m)}&=&\left(W^{(m)}\right)^t A W^{(m)} \in \R^{m\times m}\,,
\end{eqnarray*}
where $W^{(m)}$ is as defined in equation (\ref{gram}).
\begin{proposition}
\label{tri}
The matrix $T^{(m)}$ is symmetric and positive semidefinite. Furthermore
$T^{(m)}$  is tridiagonal, i.e $t_{ij}=0$ for $|i-j|\geq 2$.
\end{proposition}
\begin{proof}
The first two statements are obvious.
Let $i\leq j-2$. As $w_i \in \K^{(i)}\,$, the vector $Aw_i$ lies in the subspace $\K^{(i+1)}$. As
$j>i+1$, the vector $w_j$ is orthogonal on   $\K^{(i+1)}$, in other words
\begin{eqnarray*}
t_{ji}&=&\langle w_j,Aw_i\rangle=0\,.
\end{eqnarray*}
As $T^{(m)}$ is symmetric, we also have $t_{ij}=0$ which proves the assertion.
\end{proof}
We will see in section \ref{shrinkpls} that the matrices $T^{(m)}$ and their eigenvalues determine the shrinkage factors of
the PLS estimator. To prove this, we list some properties of $T^{(m)}$ in
teh following sections.
\section{Tridiagonal matrices}\label{sectiontri}

\begin{definition}
A symmetric tridiagonal matrix $T$ is called unreduced if all subdiagonal
entries are non-zero, i.e $t_{i,i+1} \not=0$ for all $i$.
\end{definition}
\begin{theorem}[\cite{Parlett9801}]
All eigenvalues of an unreduced matrix  are distinct.
\end{theorem}

Set
\begin{eqnarray*}
T^{(m)}&=&
\begin{pmatrix}
a_1 & b_1 & 0& \ldots & 0\\
b_1 & a_2 & b_2 &\ldots &0 \\
\ldots &\ldots &\ldots&\ldots&\vdots \\
0&0&\ldots &a_{m-1}&b_{m-1}\\
0&0 &  \ldots & b_{m-1} & a_m
\end{pmatrix} \,.
\end{eqnarray*}

\begin{proposition}\label{unreduced}
If $\dim \K^{(m)}=m$, the matrix $T^{(m)}$ is unreduced. More precisely $b_i
>0$ for all $i \in \{1,\ldots,m-1\}\,$.
\end{proposition}
\begin{proof}
Set $v_i=A^{i-1}b$ and denote  by $w_1,\ldots,w_m$
the basis obtained by Gram-Schmidt. Its existence is guaranteed as we assume
that $\dim \K^{(m)}=m$. For simplicity of notation, we assume
that the vectors $w_i$ are not normalized to have length 1. By definition
\begin{eqnarray}
\label{in2}
w_i&=&v_i -\sum_{k=1} ^{i-1} \frac{\langle v_i,w_k\rangle}{\langle
w_k,w_k\rangle}\cdot w_k\,.
\end{eqnarray}
As the vectors $w_i$ are pairwisse orthogonal,  It follows that 
\begin{eqnarray*}
\label{in1}
\langle w_i,v_i\rangle&=&\langle v_1,v_i\rangle >0\,.
\end{eqnarray*}
We conclude that
\begin{eqnarray*}
b_i&=&\langle w_i,Aw_{i-1}\rangle \\
&\stackrel{(\ref{in2})}{=}&\left \langle
w_i,A\cdot\left( v_{i-1} - \sum_{k-1} ^{i-2} \frac{\langle v_{i-1},w_k\rangle}{\langle
w_k,w_k\rangle}\cdot w_k \right)\right \rangle \\
&\stackrel{Av_{i-1}=v_i}{=}&\langle w_i,v_i\rangle - \sum_{k=1} ^{i-2} \frac{\langle v_{i-1},w_k\rangle}{\langle
w_k,w_k\rangle} \langle w_i,Aw_k \rangle\\
&\stackrel{(\ref{tri})}{=} &\langle w_i,v_i \rangle \\
&\stackrel{(\ref{in1})} = &\langle v_i,v_i \rangle >0
\end{eqnarray*}
\end{proof}
Note that the matrix
$T^{(m-1)}$ is obtained from $T^{(m)}$ by deleting the last column and row of
$T^{(m)}$. It follows that we can give a recursive formula for the
characteristical polynomials 
\begin{eqnarray*}
\chi^{(m)}:=\chi_{T^{(m)}}
\end{eqnarray*}
of $T^{(m)}$. We have
\begin{eqnarray}
\label{charpol}
\chi^{(m)}\left(\lambda\right)= \left(a_m-\lambda\right)\cdot \chi^{(m-1)}(\lambda) -
b^2 _{m-1}\chi^{(m-2)}(\lambda)
\end{eqnarray}
and $\chi^{(1)}(X)=a_1 - X\,$.\\

We want to deduce properties of the eigenvalues of $T^{(m)}$ and $A$ and explore
their relationship. Denote the eigenvalues of $T^{(m)}$ by
\begin{eqnarray}
\label{evt}
\mu^{(m)} _1 > \ldots > \mu^{(m)} _{m}\geq 0\,.
\end{eqnarray}
\begin{remark}
All eigenvalues of $T^{\left(m^*\right)}$ are eigenvalues of $A$.
\end{remark}
\begin{proof}
First note that
\begin{eqnarray*}
A_{\left| \mathcal {K}^{(m^*)}\right. }&:& \mathcal {K}^{(m^*)}
\longrightarrow \mathcal {K}^{(m^*+1)}\stackrel{\ref{dimkm} }{=}\mathcal {K}^{(m^*)}\,.
\end{eqnarray*}
As the columns of the matrix $W^{\left(m^*\right)}$ form an orthonormal basis
of $\K^{\left(m^*\right)}\,$, 
\begin{eqnarray*}
T^{(m^*)}&=& \left(W^{(m^*)}\right)^t A W^{(m^*)}
\end{eqnarray*}
is the matrix that represents $A_{\left| \mathcal {K}^{(m^*)}\right. }$ with
repect to this basis. As any eigenvalue of  $A_{\left| \mathcal
{K}^{(m^*)}\right. }$ is obviously an eigenvalue of $A$, the proof is complete
\end{proof}

The following theorem is a special form of the Cauchy Interlace Theorem.  In
this version, we use a general result from \cite{Parlett9801} and exploit the tridiagonal structure of $T^{(m)}$.

\begin{theorem}
\label{CI2}
Each interval 
\begin{eqnarray*}
\left [\mu^{(m)} _{m-j},\mu^{(m)} _{m-(j+1)}\right ]
\end{eqnarray*} $(j=0,\ldots,m-2)$ contains a different
eigenvalue of $T^{(m+k)})$ ($k\geq1$). In addition, there is a different eigenvalue of
$T^{(m+k)}$ outside the open interval $(\mu^{(m)} _{m},\mu^{(m)} _{1})\,$. 
\end{theorem}
This theorems ensures in particular that there is a different eigenvalue of $A$  in
the interval $\left [\mu^{(m)} _k, \mu^{(m)} _{k-1} \right]$. Theorem
\ref{CI2} holds independently of assumption (\ref{ass}).
\begin{proof}
By definition, for $k\geq 1$ 
\begin{eqnarray*}
T^{(m+k)} &=& \begin{pmatrix} T^{(m-1)} & t^t & 0 \\

t & a_m & * \\
0 & * & *
\end{pmatrix} \,.
\end{eqnarray*}
Here $t=(0,\ldots,0,b_{m-1})$, so 
\begin{eqnarray*}
T^{(m)} &=& \begin{pmatrix} T^{(m-1)} & t^t \\

t & a_m 
\end{pmatrix} \,.
\end{eqnarray*}
An application of Theorem 10.4.1 in \cite{Parlett9801} gives the desired result. 
\end{proof}
\begin{lemma}\label{distinct}
If $T^{(m)}$ is unreduced, the eigenvalues of $T^{(m)}$ and the
eigenvalues of $T^{(m-1)}$ are distinct.
\end{lemma}
\begin{proof}
Suppose the two matrices have a common eigenvalue $\lambda$. It follows from
(\ref{charpol}) and the fact that $T^{(m)}$ is unreduced that $\lambda$ is an eigenvalue of $T^{(m-2)}$. Repeating
this, we deduce that $a_1$ is an eigenvalue of $T^{(2)}$, a contradiction,
as
\begin{eqnarray*}
0&=&\chi^{(2)} (a_1)=-b_1 ^2 \,.
\end{eqnarray*}
\end{proof}
\begin{remark}
\label{not1}
In general it is not true  that $T^{(m)}$ and a submatrix  $T^{(k)}$ have
distinct  eigenvalues. Consider the  case where  $a_i=c$ for all $i$. Using
equation (\ref{charpol}) we conclude that 
$c$ is an eigenvalue for all submatrices with  $m$ odd.
\end{remark}
\begin{proposition}
If $\dim \K^{(m)}=m$, we have $\det \left (T^{(m-1)} \right)\not= 0$.
\end{proposition}
\begin{proof}
$T^{(m)}$ is positive semidefinite , hence all eigenvalues of $T^{(m)}$ are
$\geq0$. In other words,  $\det \left (T^{(m-1)}
\right)\not= 0$ if and only if its smallest eigenvalue $\mu^{(m-1)}_{m-1}$
is $>0$. Using Theorem \ref{CI2} we have
\begin{eqnarray*}
\mu_m ^{(m)} \geq \mu^{(m-1)}_{m-1} \geq 0\,.
\end{eqnarray*}
As $\dim \K^{(m)}=m$, the matrix $T^{(m)}$ is unreduced, which implies that $T^{(m)}$ and $T^{(m-1)}$ have no common eigenvalues (see
\ref{distinct}). We can  therefore replace the first $\geq$ by $>$, i.e. the smallest
eigenvalue of $T^{(m-1)}$ is $>0$. 
\end{proof}
In general, it is not true that $\det \left (T^{(m)} \right)\not= 0$. An easy example is 
\begin{eqnarray*}
A= \begin{pmatrix}
2 & 0\\
0 & 0
\end{pmatrix} &, &b=\begin{pmatrix} 1\\1 \end{pmatrix}\,.
\end{eqnarray*}
We have
\begin{eqnarray*}
K^{(2)}(A,b)&=&\left(b,Ab\right)\\
&=&\begin{pmatrix}
1& 2\\
1 & 0
\end{pmatrix}\,,
\end{eqnarray*}
i.e. $\dim \K^{(2)}=2$. On the other hand 
\begin{eqnarray*}
\det \left(T^{(2)}\right) &=& \det \begin{pmatrix} 1 & 1 \\ 1 &1 \end{pmatrix} =0\,.
\end{eqnarray*}

It is well known that the matrices $T^{(m)}$ are closely related to the
so-called Rayleigh-Ritz procedure, a method that is used to approximate
eigenvalues. For details consult e.g. \cite{Parlett9801}. 

\section{What is shrinkage?}
\label{shrinkage}
We have presented two estimators for the regression parameter $\beta$ -- OLS
and PLS -- which
also define estimators for $X\beta$ via 
\begin{eqnarray*}
\hat y_{\bullet}&=&X\cdot \hat \beta_{\bullet}\,.
\end{eqnarray*}
One possibility to evaluate the quality of an estimator is to determine
its  Mean Squared Error (MSE). In general, the MSE of an
estimator $\hat \theta$ for a vector-valued parameter $\theta$ is defined
as
\begin{eqnarray*}
\text{MSE}\left(\hat\theta\right)&=&E\left[\text{trace}\left(\hat\theta
-\theta\right)\left(\hat \theta -\theta\right)^t\right] \\
&=& E\left[ \left(\hat\theta
-\theta\right)^t\left(\hat \theta -\theta\right) \right]\\
&=& \left(E\left[ \hat \theta \right] -\theta \right)^t\left(E\left[ \hat
\theta \right] -\theta \right)+ E\left[ \left(\hat \theta^t - E\left[\hat\theta\right]\right)^t
\left(\hat \theta^t - E\left[\hat\theta\right]\right) \right]\,.
\end{eqnarray*}
This is the well-known bias-variance decomposition of the MSE. The first
part is the squared bias  and the second part is the variance term.

We start by investigating the  class of linear estimators, i.e.
estimators that are of the form $\hat \theta= Sy$ for some matrix $S$ that
does not depend on $y$. The OLS estimators are linear:
\begin{eqnarray*}
\hat \beta_{OLS}&=&\left(X^t X \right)^- X^t y:=S_1 y\\
\hat y_{OLS}&=& X\cdot\left(X^t X \right)^- X^t  y:=S_2 y  \,.
\end{eqnarray*}
$S_2$ is the projection $P_{L(X)}$ onto the space that is spanned by the columns of
$X$.\\

Recall the regression model (\ref{linreg}). 
\begin{proposition}
Let $\hat \theta=Sy$ be a linear estimator. We have
\begin{eqnarray*}
E\left[\hat \theta \right]&=& SX\beta \\
\text{var}\left[\hat \theta \right]&=&\sigma^2 \text{tr}\left(SS^t\right)\,.
\end{eqnarray*}
\end{proposition}
The estimator $\hat y_{OLS}$ is
unbiased as 
\begin{eqnarray*}
E\left[\hat y_{OLS}\right] &=& S_2 X\beta \\
&=& P_{L(X)} X\beta\\
&=& X\beta\,.
\end{eqnarray*}
The estimator $\hat \beta_{OLS}$ is only unbiased if  $\beta \in \text{range}
\left(X^tX\right)^- \,$:
\begin{eqnarray*}
E\left[\hat\beta_{OLS}\right]&=&E\left[\left(X^tX\right)^- X^t y\right]\\
&=& \left(X^tX\right)^- X^t E\left[y\right]\\
&=&  \left(X^tX\right)^- X^t X\beta \\
&=& \beta\,.
\end{eqnarray*}

Let us now have a closer look at the variance term. 

For $\hat \beta_{OLS}$ we have
\begin{eqnarray*}
S_1S_1 ^t&=& \left(X^tX\right)^- X^t X \left(X^t X\right)^- \\
&=& \left(X^tX\right)^-\\
&=& U \Lambda^- U^t \,,
\end{eqnarray*}
hence 
\begin{eqnarray}
\label{varbeta}
var\left(\hat \beta_{OLS}\right) &=& \sigma^2 \cdot  \sum_{i=1} ^{p^*} \frac{1}{\lambda_i}\,.
\end{eqnarray}
Next note that $S_2$ is the operator that projects on the space spanned by the columns of
$X$. It follows that $\text{tr}(S_2S_2 ^t)=\text{rk}(X)=p^*$ and that
\begin{eqnarray*}
var\left(\hat y_{OLS}\right) &=& \sigma^2 \cdot p^*\,.
\end{eqnarray*}
We conclude that the MSE of the estimator $\hat \beta_{OLS}$ depends on the
eigenvalues $\lambda_1,\ldots,\lambda_{p^*}$ of  $A=X^tX$. Small eigenvalues of $A$
correspond to directions in $X$ that have very low variance. Equation
(\ref{varbeta})  shows that if some eigenvalues are small, the variance of
$\hat \beta_{OLS}$ is very high, which leads to a high MSE. \\

One possibility  to (hopefully) decrease the MSE is to modify the OLS
estimator by shrinking the directions of the OLS estimator
that are responsible for a high variance. This of course introduces bias.
We shrink the OLS estimator in the hope that the increase in bias is small
compared to the decrease in variance. \\

In general, a  shrinkage estimator for $\beta$ is of  the form
\begin{eqnarray*}
\hat \beta_{shr} &= &\sum_{i=1} ^{p^*} f(\lambda_i) \cdot z_i\,,
\end{eqnarray*}
where $f$ is some real-valued function.  The values $f(\lambda_i)$ are
called shrinkage factors.

Examples are
\begin{itemize}
\item Principal Component Regression
\begin{eqnarray*}
f(\lambda_i)&=&\begin{cases}
1 & \text{ith principal component included} \\
0 & \text{otherwise}
\end{cases}
\end{eqnarray*}
and
\item Ridge Regression
\begin{eqnarray*}
f(\lambda_i)&=&\frac{\lambda_i}{\lambda_i +\lambda}
\end{eqnarray*}
where $\lambda>0$ is the Ridge parameter.
\end{itemize}
We will see in section \ref{shrinkpls} that PLS is a shrinkage estimator as
well. It will turn out that the shrinkage behavior of PLS regression is
rather complicated. \\

Let us investigate in which way the MSE of the estimator is influenced by
the shrinkage factors. If the shrinkage estimators are linear, i.e. the
shrinkage factors do not depend on $y$, this is an easy task. Let us first
write the shrinkage estimator in matrix notation. We have 
\begin{eqnarray*}
\hat \beta_{shr} &=& S_{shr,1} y\\
&=& U \Sigma^- D_{shr} V^t y \,.
\end{eqnarray*}
The diagonal matrix $D_{shr}$ has entries
$f(\lambda_i)$. The shrinkage estimator for $y$ is 
\begin{eqnarray*}
\hat y_{shr} &=& S_{shr,2} y\\
&=& V\Sigma \Sigma^- D_{shr} V^t\,.
\end{eqnarray*}
We calculate the variance of these estimators. 
\begin{eqnarray*}
\text{tr}\left(S_{shr,1} S_{shr,1} ^t\right)&=& \text{tr}\left(  U
\Sigma^-D_{shr} \Sigma^-D_{shr} U^t\right)\\
&=& \text{trace}\left(\Sigma^-D_f \Sigma^-D_f  \right)\\
&=& \sum_{i=1} ^{p^*} \frac{\left(f(\left(\lambda_i\right)\right)^2}{\lambda_i}
\end{eqnarray*}
and
\begin{eqnarray*}
\text{tr}\left(S_{shr,2} S_{shr,2} ^t\right)&=&\text{tr}\left(V\Sigma
\Sigma^- D_{shr} \Sigma \Sigma^- D_{shr} V^t\right)\\
&=& \text{tr}\left(\Sigma
\Sigma^- D_{shr} \Sigma \Sigma^- D_{shr} \right)\\
&=& \sum_{i=1} ^{p^*} \left(f(\left(\lambda_i\right)\right)^2\,.
\end{eqnarray*}
Next, we calculate the bias of the two shrinkage estimators.
 We have 
\begin{eqnarray*}
E\left[S_{shr,1} y\right]&=& S_{shr,1} X\beta\\
&=& U\Sigma D_{shr}  \Sigma^-U^t \beta
\,.
\end{eqnarray*}
It follows that 
\begin{eqnarray*}
\text{bias}^2\left( \hat \beta_{shr} \right)&=& \left(E\left[S_{shr,1} y\right] -\beta
\right)^t\left(E\left[S_{shr,1} y\right] -\beta   \right)\\
&=& \left(U^t \beta\right)^t \left(\Sigma D_f  \Sigma^- -\text{Id}\right)^t  \left(\Sigma D_f  \Sigma^- -\text{Id}\right)\left(U^t \beta\right)\\
& =&\sum_{i=1} ^{p^*} \left(f(\lambda_i)
-1\right)^2 \left(u_i ^t \beta \right)^2\,.
\end{eqnarray*}
Replacing $S_{shr,1}$ by  $S_{shr,2}$ it is as easy to show that
\begin{eqnarray*}
\text{bias}^2\left( \hat y_{shr} \right)&=& \sum_{i=1} ^p \lambda_i \left(f(\lambda_i)
-1\right)^2 \left(u_i ^t \beta \right)^2\,.
\end{eqnarray*}
\begin{theorem}
For the shrinkge estimator $\hat \beta_{shr}$ and $\hat y_{shr}$ defined above we
have
\begin{eqnarray*}
MSE\left(\hat \beta_{shr} \right)&=& \sum_{i=1} ^{p^*} \left(f(\lambda_i)
-1\right)^2 \left(u_i ^t \beta \right)^2 + \sigma^2  \sum_{i=1} ^{p^*}
\frac{\left(f\left(\lambda_i\right)\right)^2}{\lambda_i} \\
MSE\left(\hat y_{shr} \right)&=& \sum_{i=1} ^{p^*} \lambda_i \left(f(\lambda_i)
-1\right)^2 \left(u_i ^t \beta \right)^2 + \sigma^2 \sum_{i=1} ^{p^*} \left(f\left(\lambda_i\right)\right)^2\,.
\end{eqnarray*}
\end{theorem}

If the shrinkage factors are deterministic, i.e. they do not depend on $y$,  any value $f(\lambda_i)\not=1$
increases the bias. Values $\left|f(\lambda_i)\right| <1$ decrease the
variance, whereas values   $\left|f(\lambda_i)\right| >1$ increase the
variance. Hence an absolute value $>1$ is always undesirable. The situation is completely different for stochastic
shrinkage factors. We will discuss this in the following section.\\

Note that there is a different notion of shrinkage, namely that the $l_2$-
norm of an estimator  is smaller than the $l_2$-norm of the OLS
estimator. Why is this a desirable property? Let us again consider the case
of linear estimators. Set $\hat \beta_i =S_i y$ for $i=1,2$. We have
\begin{eqnarray*}
\left\|\beta_i \right\|_2 ^2&=& y^t S_i ^t S_i y\,.
\end{eqnarray*}
The property that for all $y \in \mathbb{R}^n$ 
\begin{eqnarray*}
\left\|\beta_1 \right\|_2 &\leq& \left\|\beta_2 \right\|_2
\end{eqnarray*}
is  equivalent  to the
condition that
\begin{eqnarray*}
S_1 ^t S_1 -S_2 ^t S_2
\end{eqnarray*}
is negative semidefinite. The trace of negative semidefinite matrices is $\leq 0$. Furthermore
$\text{trace}\left(S_i ^t S_i\right)=\text{trace}\left(S_i  S_i^t \right)$,
so we  conclude that 
\begin{eqnarray*}
\text{var}\left(\hat\beta_1\right)&\leq&\text{var}\left(\hat\beta_2\right)\,.
\end{eqnarray*}
It is known (see \cite{Goutis9601}) that
\begin{eqnarray*}
\|\hat \beta^{(1)} _{PLS}\|_1\leq \|\hat \beta^{(2)} _{PLS}\|_2 \leq \ldots
\leq \|\hat \beta^{(m^*)}_{PLS}\|_2=\|\hat \beta_{OLS}\|_2\,.
\end{eqnarray*}

\section{The shrinkage factors of PLS}\label{shrinkpls}
In this section, we give a simpler and clearer proof of the shape of the
shrinkage factors of PLS. Basically, we combine the results of \cite{Butler0001}
and \cite{Phatak0301}. It turns out that some of the factors
$f^{(m)}(\lambda_i)$ are greater than 1. We try to explain why these "peculiar
shrinkage properties" do not necessarily imply that the MSE of the PLS
estimator is increased.\\

Denote by $\pi^{(m)}$ the polynomial associated to $T^{(m)}$ that was
defined in proposition \ref{penrose}, i.e. 
\begin{eqnarray*}
\pi^{(m)}\left(T^{(m)}\right)=\pi_{T^{(m)}} \left(T^{(m)}\right)= \left(T^{(m)}\right)^-\,.
\end{eqnarray*}
Recall that the eigenvalues of $T^{(m)}$ are denoted by $\mu_m ^{(m)}$. It
follows that
\begin{eqnarray}
\label{TM}
f^{(m)}(\lambda)&:=&\lambda\cdot \pi^{(m)}(\lambda)=1-\prod_{i=1} ^m \left(1-\frac{\lambda}{\mu^{(m)} _i} \right)\,.
\end{eqnarray}
By definition of PLS, $\hat \beta_{PLS} ^{(m)} \in \K^{(m)}$ hence
there is a polynomial $\pi$ of degree $\leq m-1$ with $\hat \beta_{PLS}
^{(m)}=\pi(A)b$.
\begin{proposition}[\cite{Phatak0301}]
Suppose that $\dim \K^{(m)}=m$. We have
\begin{eqnarray*}
\hat \beta_{PLS} ^{(m)}&=&\pi^{(m)}(A)\cdot b\,.
\end{eqnarray*}
\end{proposition}
\begin{proof}[Proof (\cite{Phatak0301})]
By proposition \ref{penrose}, 
\begin{eqnarray*}
\left(T^{(m)}\right) ^{-}&=&\pi^{(m)}\left(T^{(m)}\right)\,.
\end{eqnarray*}
We plug this into
equation (\ref{bhat}) and obtain
\begin{eqnarray*}
\hat \beta^{(m)} _{PLS} &=& W^{(m)} \pi^{(m)} \left ( \left
(W^{(m)}\right)^t A W^{(m)}  \right) \left (W^{(m)}\right)^t b\,.
\end{eqnarray*}
Recall that the columns of $W^{(m)}$ form an orthonormal basis of
$\mathcal{K}^{(m)}(A,b)$. It follows that $W^{(m)} \left(W^{(m)}\right)^t$ is
the operator that projects on the space $\mathcal{K}^{(m)}(A,b)$. In particular
\begin{eqnarray*}
W^{(m)} \left(W^{(m)}\right)^t A^j b= A^j b
\end{eqnarray*}
for $j=1,\ldots,m-1$. This implies that
\begin{eqnarray*}
\hat \beta^{(m)} _{PLS} &=& \pi^{(m)}(A)\cdot b\,.
\end{eqnarray*}
\end{proof}
\begin{corollary}[\cite{Phatak0301}]Suppose that $\dim \K^{(m)}=m$.
If we denote by $z_i$ the component of $\hat{\beta}_{OLS}$ along the
$i$th eigenvector of $A$ then
\begin{eqnarray*}
\hat{\beta}_{PLS} ^{(m)} &=& \sum_{i=1} ^{p^*} f^{(m)}(\lambda_i) \cdot z_i\,,
\end{eqnarray*}
where $f^{(m)}$ is the polynomial  defined in (\ref{TM}).
\end{corollary}
\begin{proof}(\cite{Phatak0301})
This follows immediately from the proposition above. We have
\begin{eqnarray*}
\hat \beta^{(m)} _{PLS} &=& \pi^{(m)}(A)b \\
&=& U \pi^{(m)} (\Lambda) \Sigma V^t y\\
&=& \sum_{i=1} ^{p^*} \pi^{(m)}(\lambda_i) \sqrt{\lambda_i}(v_i)^t y u_i\\
&=& \sum_{i=1} ^p \pi^{(m)}(\lambda_i) \lambda_i
\frac{1}{\sqrt{\lambda_i}}(v_i ^t y) u_i\\
&\stackrel{(\ref{TM})}{=}& \sum_{i=1} ^p f^{(m)}(\lambda_i) z_i\,.
\end{eqnarray*}
\end{proof}
We now show that some of the shrinkage factors of PLS are $\not=1\,$. 
\begin{theorem}[\cite{Butler0001}]
For each $m\leq m^* -1$, we can decompose the interval
$\left[\lambda_p,\lambda_i\right]$ into $m+1$ disjoint intervals\footnote{We say that $I_j \leq  I_k$ if $\sup I_j \leq \inf I_k\,$.} 
\begin{eqnarray*}
I_1 \leq I_2 \leq \ldots \leq I_{m+1}
\end{eqnarray*}
 such that 
\begin{eqnarray*}
f^{(m)} \left(\lambda_i \right) \begin{cases}
 \leq 1 &\lambda_i \in I_j \text { and } j \text{ odd} \\
\geq 1 & \lambda_i \in I_j \text{ and } j \text{ even} 
\end{cases}
\,.
\end{eqnarray*}

\end{theorem}

\begin{proof}
Set $g^{(m)}=1-f^{(m)}$. It follows from equation (\ref{TM}) that the zero's
of $g^{(m)}$ are  $\mu^{(m)} _m,\ldots,\mu^{(m)} _1$. As $T^{(m)}$ is
unreduced, all eigenvalues are  distinct. Set $\mu^{(m)}_{0}=\lambda_1$
and $\mu^{(m)}_{m+1}=\lambda_p$.  Define $I_j=
]\mu_{i} ^{(m)},\mu_{i+1} ^{(m)}[$ for $j=0,\ldots,m$\,.  By definition,
$g^{(m)}(0)=1$. Hence $g^{(m)}$ is
non-negative on the intervals $I_j$ if $j$ is odd and $g^{(m)}$ is non-positive  on the
intervals $I_j$  if $j$ is even. It follows
from Theorem \ref{CI2}  that all interval $I_j$ contain at least one
eigenvalue $\lambda_i$ of $A\,$.
 \end{proof}
In general it is not true that   $f^{(m)}(\lambda_i)\not=1$ for all
$\lambda_i$ and $m=1,\ldots,m^*\,$.  Using the example in remark \ref{not1}
and the fact that
\begin{eqnarray*}
f^{(m)}(\lambda_i)&=&1
\end{eqnarray*}
is equivalent to the condition that $\lambda_i$ is an eigenvalue of
$T^{(m)}$, it is  easy to  construct a counterexample. Using some of the results of section \ref{sectiontri}, we can however deduce
that some factors are indeed $\not=1$.  As all eigenvalues of $T^{(m^*-1)}$
and $T^{(m^*)}$ are distinct (c.f. proposition \ref{distinct}), we see that
$f^{(m^*-1)}(\lambda_i) \not=1$ for all $i$. In particular 
\begin{eqnarray*}
f^{(m^*-1)}(\lambda_1) 
\begin{cases}
<1 & m^* \text{ even }\\
>1 &   m^* \text{ odd}
\end{cases}\,.
\end{eqnarray*}
More generally, using proposition \ref{distinct}, we conclude that
$f^{(m-1)}\left(\lambda_i\right)$ and  $f^{(m)}\left(\lambda_i\right)$ is
not possible. In practice -- i.e. calculated on a data set --  the factors seem to be $\not=1$ all of the time. \\
Furthermore
\begin{eqnarray*}
0\leq f^{(m)}(\lambda_p)<1\,.
\end{eqnarray*}
To proove this, we set $g^{(m)}=1-f^{(m)}$. We have by definition $g^{(m)}(0)=1$. Furthermore, the smallest
positive zero of $g^{(m)}$ is $\mu_m ^{(m)}$ and it follows from  Theorem
\ref{CI2} and proposition \ref{distinct} that $\lambda_p < \mu_m ^{(m)}$. Hence $g^{(m)}(\lambda_p)
\in ]0,1]$.\\

Using Theorem \ref{CI2}, more precisely
\begin{eqnarray*}
\lambda_p \leq \mu_i ^{(m)} \leq \lambda_i
\end{eqnarray*}
it is possible to bound the terms
\begin{eqnarray*}
1-\frac{\lambda_i}{\mu_i ^{(m)}}\,.
\end{eqnarray*}
From this we can derive bounds on the shrinkage factors. We will not pursue
this further, readers who are interested in the bounds should consult
\cite{Lingjaerde0001}. Instead, we have a closer look at the MSE of the PLS estimator.\\

In section \ref{shrinkage} we showed that a value $|f^{(m)}(\lambda_i)|>1$ is not
desirable, as the variance of the estimator increases.  Note however, that in
the case of PLS, the factors $f^{(m)}(\lambda_i)$ are stochastic; they depend on
$y$ - in a nonlinear way. For $\hat \beta_{PLS} ^{(m)}$ we have the following situation:  If we set
$Z=f^{(m)}(\lambda_i)$ and $W= \frac{\left(v_i\right)^ty}{\sqrt{\lambda_i}}$, we have to compare
\begin{eqnarray*}
var(Z\cdot W) &\,\,\,\,\,\, \text{to}& var(W)\,.
\end{eqnarray*}
Note that the RHS is not necessarily smaller than the LHS, even if
$P(Z>1)=1$. An easy counterexample is $Z=\frac{1}{W}\,$ -- the LHS is $0$.\\

Among others, \cite{Frank9301} proposed to bound  the shrinkage factors of
the PLS estimator in
the following way. Set
\begin{eqnarray*}
\tilde f ^{(m)}(\lambda_i)&=&\begin{cases}
+1&f^{(m)}(\lambda_i)>+1\\
-1& f^{(m)}(\lambda_i)<-1\\
f^{(m)}(\lambda_i)&\text{otherwise}
\end{cases}
\end{eqnarray*}
and define a new estimator:
\begin{eqnarray}
\label{bound}
\hat \beta ^{(m)} _{BOUND}&:=& \sum_{i=1} ^p \tilde f
^{(m)}(\lambda_i) z_i \,.
\end{eqnarray}
{\par\centering\resizebox*{9cm}{6cm}{\rotatebox{270}{\includegraphics{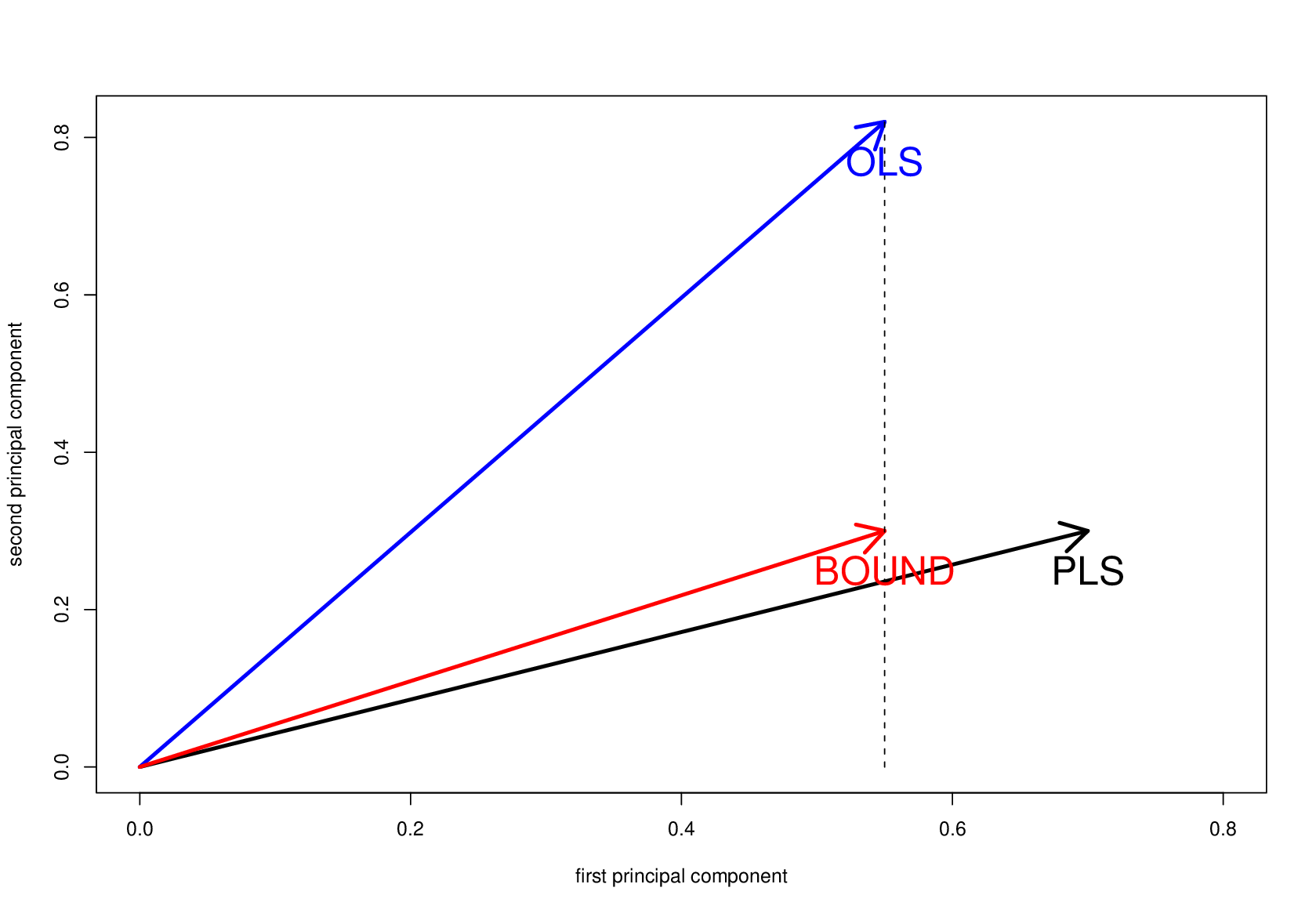}}}\par}

If the shrinkage factors are numbers, this will improve the MSE (cf. section
\ref{shrinkage}). But in the
case of stochastic shrinkage factors, the situation is completely unclear. Consider
again the example $Z=\frac{1}{W}$. Set 
\begin{eqnarray*}
\tilde Z &=&\begin{cases}
+1&Z>1\\
-1& Z<-1\\
Z&\text{otherwise}
\end{cases}
\end{eqnarray*}
In this case
\begin{eqnarray*}
0=var(Z \cdot W) < var(\tilde Z \cdot W)
\end{eqnarray*}

so  it is not clear
whether the modified estimator BOUND leads to a lower  MSE,  which was conjectured in e.g. \cite{Frank9301}.\\

The above example (involving $W$ and $Z$) is of course purely artificial. It is not  clear whether
the shrinkage factors behave this way. It is hard if not infeasable to
derive statistical properties of the PLS estimator or its shrinkage factors,
as they depend on $y$ in a complicated, nonlinear way. As an alternative, we
compare the two different estimators on different data. 
\section{Experiments}\label{Simulation}
In this section, we explore the difference between the methods  PLS
and BOUND. We investigate  three artificial datasets and one real world
example. In all examples, we rescale $X$ and $y$ to have zero mean and unit variance.\\

Let us start with the artificial datasets. Of course, artificial datasets
do not reflect many real world situations, but we have the advantage that  we
know the true regression coefficient $\beta$ and that we have an unlimited amount
of examples at hand. We can estimate the MSE of any of the
four estimators: For $k=1,\ldots,K$ we generate a sample $y$ and calculate
the estimator $\hat \theta _k$. We define
\begin{eqnarray*}
\widehat {MSE}(\hat \theta)&=& \frac{1}{K} \sum_{k=1} ^K \left(\hat \theta_k
-\theta\right)^t  \left(\hat \theta_k
-\theta\right)\,.
\end{eqnarray*}
For all examples, we choose $K=200\,$.
\subsection*{First example}
In our first example we generate $n=30$ examples in the following way: The
input data is the realistion of a $p=10$ dimensional normally distributed
variable with expectation $\bf{0}\in \mathbb{R}^p$ and covariance matrix $\Sigma
\in \mathbb{R}^{p\times p}$
defined as 
\begin{eqnarray*}
\Sigma_{ij} =\begin{cases} 1.5 & i=j\\
1 &i\not=j \end{cases}\,.
\end{eqnarray*}
 The regression coefficient $\beta$ is the random permutation of $(0,0,0,0,0,z_1,\ldots,z_5)$ with $z_i \sim
N(2,2^2)$. \\

Next we determine the variance of the error term.  We do this by considering several signal-to-noise-ratios (stnr). This
quantity is defined as
\begin{eqnarray*}
stnr&=&\frac{var(X\beta)}{var(\varepsilon)}\,.
\end{eqnarray*}
We set $stnr=1,4,16$ and determine the corresponding value of
$\sigma\,$.  We generate $K=200$ samples $y$ and calculate the four
estimators. \\

The following figures show the estimated MSE for $\beta$ and $X\beta$
respectively. The solid lines with the $\bullet$'s correspond to PLS. the lines
with the $+$'s correspond to BOUND.\\
\newpage
\begin{figure}[h]
{\par\centering\resizebox*{12cm}{3.7cm}{\rotatebox{270}{\includegraphics{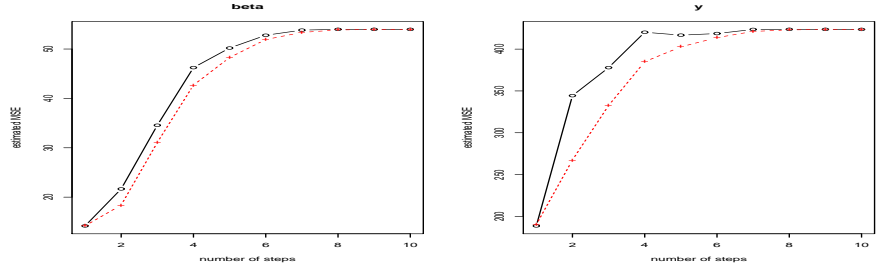}}}\par}
\caption{First example: Comparison of PLS and BOUND ($stnr=1$)}
\end{figure}
\begin{figure}[h]
{\par\centering\resizebox*{12cm}{3.7cm}{\rotatebox{270}{\includegraphics{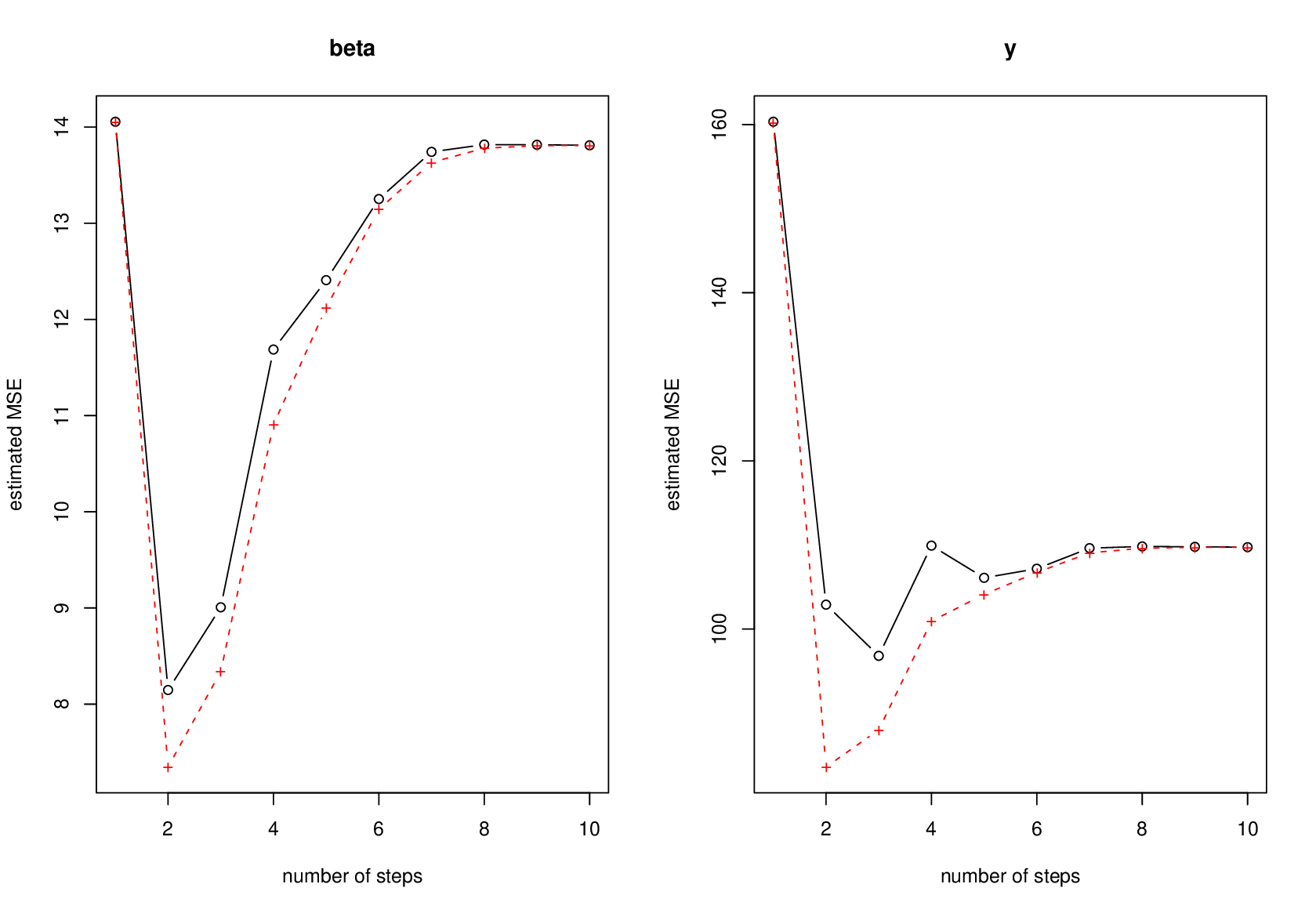}}}\par}
\caption{First example: Comparison of PLS and BOUND ($stnr=4$)}
\end{figure}

\begin{figure}[h]
{\par\centering\resizebox*{12cm}{3.7cm}{\rotatebox{270}{\includegraphics{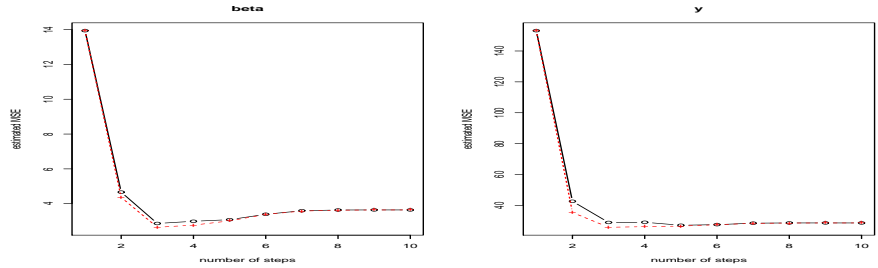}}}\par}
\caption{First example: Comparison of PLS and BOUND ($stnr=16$)}
\end{figure}

We see that BOUND is better in all cases, although the improvement is not
dramatic. We should remark that both method pick the same (optimal) number
of steps most of the times. The difference between the two methods is
especially tiny (but non-zero)
 in the first step. We do not have an explanation for this phenomenon. The
MSE is the same for the last step $m=10$ as in this case
\begin{eqnarray*}
\hat \beta_{PLS} ^{(m)} = \hat \beta_{BOUND} ^{(m)} = \hat \beta_{OLS}\,.
\end{eqnarray*}
\subsection*{Second example}
In this example, we generate $n=40$ examples. The input data is the
realisation of a $p=20$ dimensional random variable with distribution
$N(0,\Sigma)$. The covariance matrix is defined as in the first example
(with $p=10$ replaced by $p=20$). Again, the coefficients of $\beta$ are a
random permutation $\left(0,\ldots,0,z_1,\ldots,z_10\right)$ with  $z_i \sim
N(2,2^2)\,$. We consider
the  signal-to-noise-ratios $1,4,16\,$.\\
\newpage

\begin{figure}[htb]
{\par\centering\resizebox*{12cm}{3.7cm}{\rotatebox{270}{\includegraphics{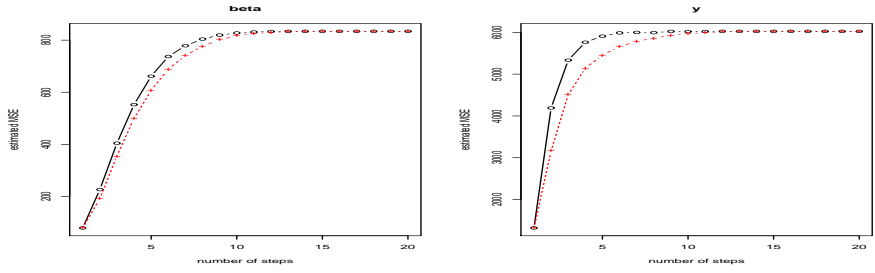}}}\par}
\caption{Second example: Comparison of PLS and BOUND ($stnr=1$)}
\end{figure}
\begin{figure}[htb]
{\par\centering\resizebox*{12cm}{3.7cm}{\rotatebox{270}{\includegraphics{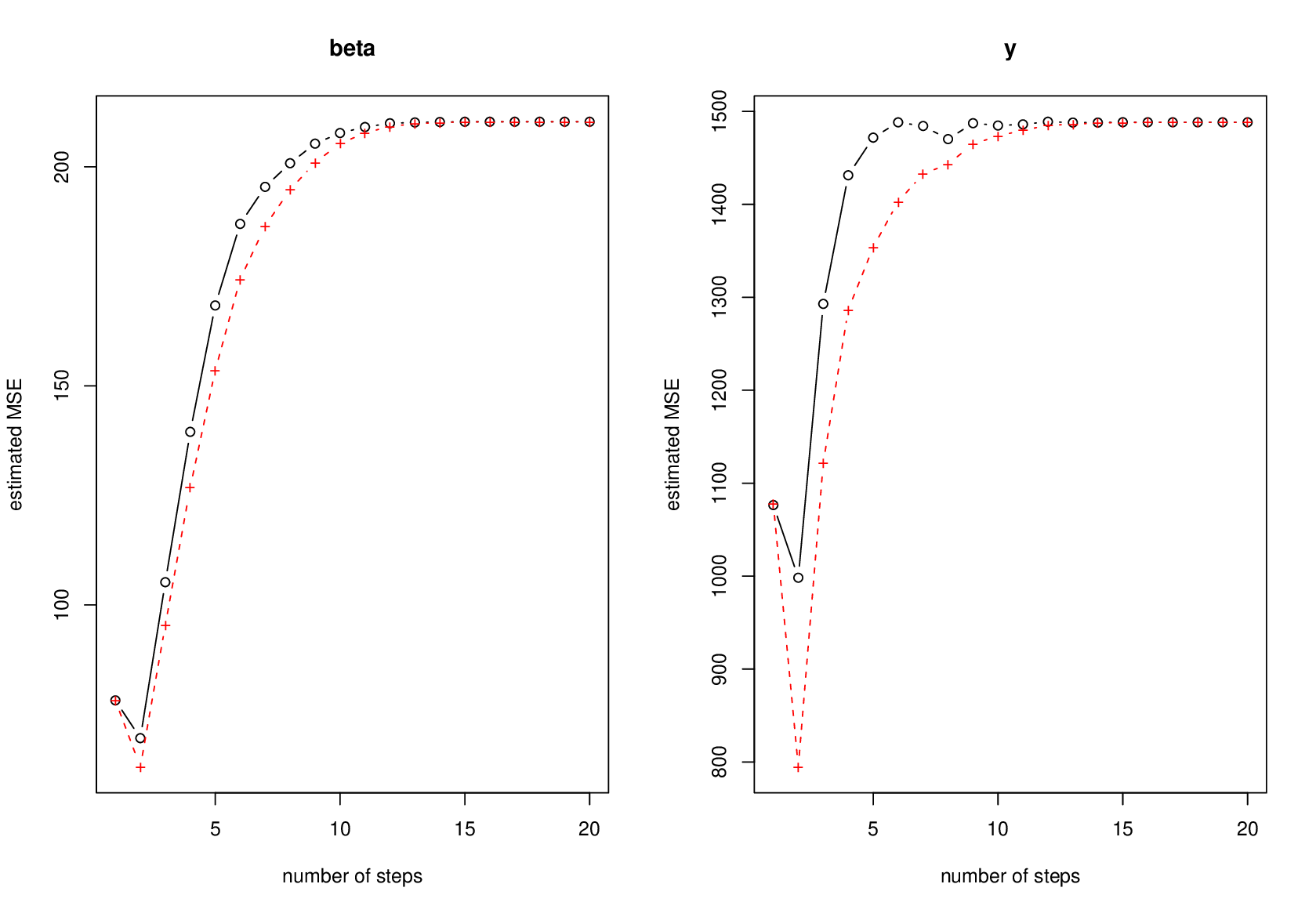}}}\par}
\caption{Second example: Comparison of PLS and BOUND ($stnr=4$)}
\end{figure}
\begin{figure}[htb]
{\par\centering\resizebox*{12cm}{3.7cm}{\rotatebox{270}{\includegraphics{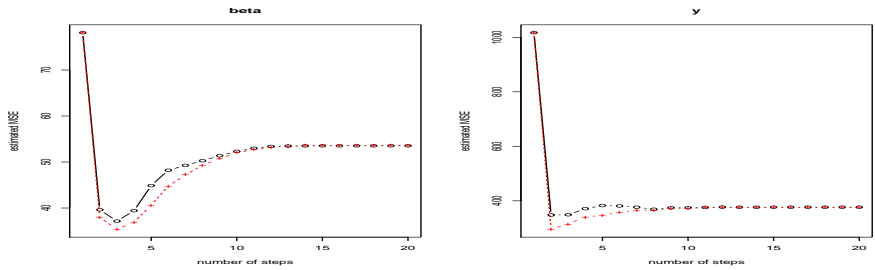}}}\par}
\caption{Second example: Comparison of PLS and BOUND ($stnr=16$)}
\end{figure}
The results are qualitatively the same as those from the first
example. BOUND is better all of the times, the optimal number of steps are
the same for both methods.
\subsection*{Third example}
The input data is generated as in the second example, in particular, we have
$p=20\,$.  This time, we only
generate $n=10$ examples. The coefficients of the regression vector $\beta$
are realizations of a $N(2,2^2)$ distibuted random variable. We investigate
the signal-to-noise-ratios $1,4,16\,$. As we have more
variables than examples, we do not  investigate estimators for
$\beta\,$: Different vectors $\beta_1 \not= \beta_2$ can lead to
$X\beta_1=X\beta_2$, so  it does not make sense to determine the bias of an
estimator for $\beta\,$. Instead, we only show the figures for $\hat y_{PLS}$ and
$\hat y_{BOUND}$. 
\newpage

\begin{figure}[htb]

{\par\centering\resizebox*{10cm}{3.7cm}{\rotatebox{270}{\includegraphics{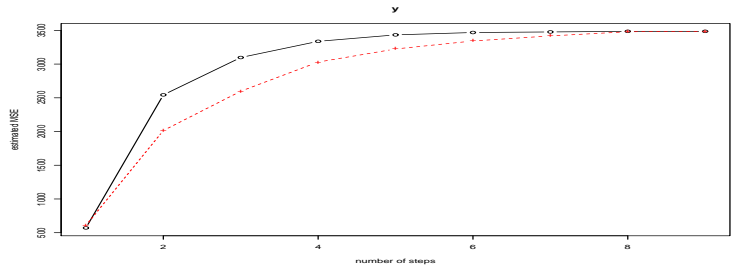}}}\par}
\caption{Third example: Comparison of PLS and BOUND ($stnr=1$)}
\end{figure}
\begin{figure}[htb]

{\par\centering\resizebox*{10cm}{3.7cm}{\rotatebox{270}{\includegraphics{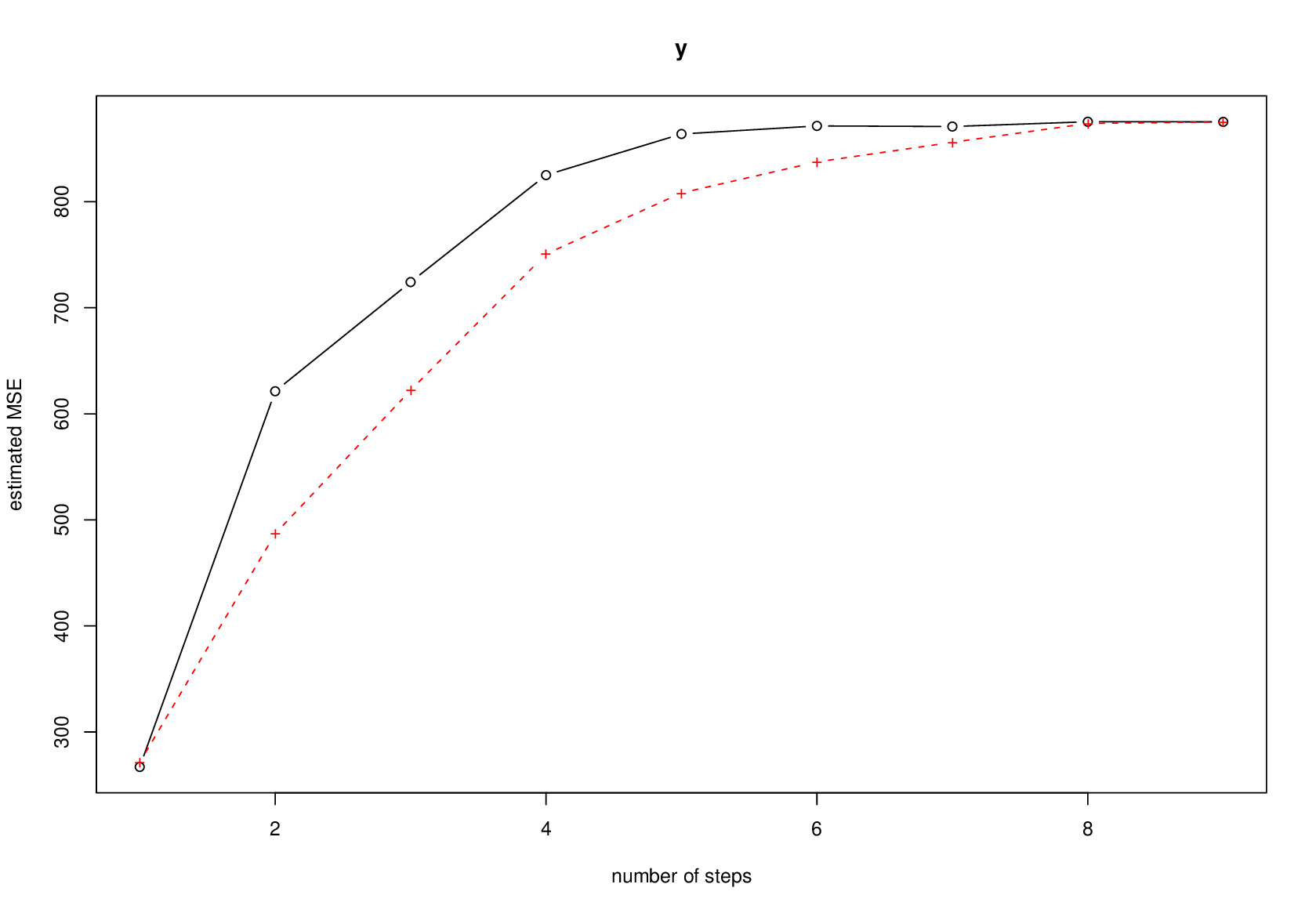}}}\par}
\caption{Third example: Comparison of PLS and BOUND ($stnr=4$)}
\end{figure}
\begin{figure}[htb]

{\par\centering\resizebox*{10cm}{3.7cm}{\rotatebox{270}{\includegraphics{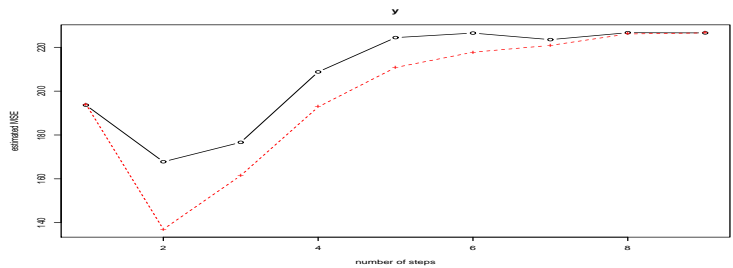}}}\par}
\caption{Third example: Comparison of PLS and BOUND ($stnr=16$)}
\end{figure}
Again, the estimated MSE of BOUND is lower than the estimated MSE of PLS.
\subsection*{Fourth example}
This example is taken from \cite{Martens0101}. A survey investigated the
degree of job satisfaction of the employees of a company. The
employees filled in a questionnaire that consisted of $p=26$ questions
regarding their work environment and one question (the response variable) regarding the degree to
which they are satisfied with their job. The answers of the employees were
summerized for each of the $n=34$ departments of the company. \\

We compare the two methods PLS and BOUND on this data set. For each
$m=1,\ldots 26$ we determine the 10fold crossvalidation error.
\newpage
\begin{figure}[htb]
\label{cv}
{\par\centering\resizebox*{12cm}{5cm}{\rotatebox{270}{\includegraphics{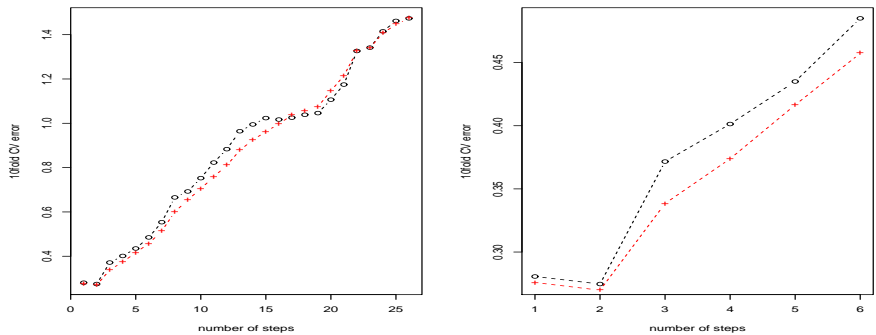}}}\par}
\caption{Left: 10fold crossvalidation error. Right: 10fold crossvalidation
error for the first 6 components}
\end{figure}
The method BOUND is slightly better than PLS on this data set: The cv error
for the optimal number of components (which is $m_{opt}=2$) is  0.2698  for BOUND
and 0.2747 for PLS. It is remarkable that in this example the cv error of
BOUND exceeds the cv error of PLS in some cases. It is not clear if this is
due to the small number of examples (which makes the estimation unprecise)
or if this can also happen "in theory".


\section{Conclusion}
This paper consists of two parts. In the first part, we gave alternative and
hopefully clearer proofs of the
shrinkage factors of PLS. In particular, we derived the fact that some of
the shrinakge factors are $>1$. We explained in detail that this would lead
to an unnecessarily high MSE if PLS was a linear estimator. This is however
not the case and we emphasized that bounding the absolute value of the
shrinkage factors by $1$ does not automatically lead to a lower MSE. \\

In the second part, we investigated the problem numerically.  Experiments on
simulated and real world data showed that it might be  better to adjust the shrinkage factors so that their
absolute value is $\leq 1$ - a method that we called BOUND. The difference
between BOUND and PLS was not dramatic however. Besides, the scale of the
experiments was of course way too small, so it would be light-headed if we
concluded  that we should always use BOUND  instead  of PLS. \\

Nevertheless, the experiments show that it is worth exploring the method
BOUND in more detail. One drawback of this method is that we have to adjust
the shrinkage factors "by hand". If bounding the shrinkage factors tends to lead
to better results, we might modify the original optimization problem of PLS such
that the shrinkage factors of the solution are bounded.  We might modfify $A$ and $b$
to obtain a different Krylov space or replace $\K^{(m)}$ by a different set
of feasible solutions.  
 \subsection*{Acknowledgement}
I would like to thank Ulrich Kockelkorn who eliminated innumerable  errors from
earlier versions of this paper and who gave a lot  of helpful remarks. I would
also like to thank J\"org Betzin for our extensive discussions on  PLS.

\bibliographystyle{plain}

\end{document}